\documentclass[a4paper]{amsart}

\usepackage{amsmath}
\usepackage{amssymb}

\usepackage[all]{xy}
\def\mod{\operatorname{mod}}
\def\add{\operatorname{add}}
\def\Hom{\operatorname{Hom}}
\def\Ext{\operatorname{Ext}}
\def\End{\operatorname{End}}
\def\dim{\operatorname{dim}}
\def\CC{\mathcal{C}}
\def\DD{\mathcal{D}}
\def\SS{\mathcal{S}}
\def\ol{\overline}
\def\op{\textrm{op}}

\newtheorem{thm}{Theorem}[section]
\newtheorem{prop}[thm]{Proposition}
\newtheorem{lem}[thm]{Lemma}
\newtheorem{cor}[thm]{Corollary}
\newtheorem*{Thm}{Theorem}

\title[Derived equivalence for Cluster-tilted algebras of type
$A_n$]{Derived equivalence classification for Cluster-tilted algebras
  of type $A_n$}

\author{Aslak Bakke Buan}
\address{Institutt for matematiske fag\\
Norges teknisk-naturvitenskapelige universitet\\
N-7491 Trondheim\\
Norway}
\email{aslakb@math.ntnu.no}

\author{Dagfinn F. Vatne}
\address{Instututt for matematiske fag\\
Norges teknisk-naturvitenskapelige universitet\\
N-7491 Trondheim\\
Norway}
\email{dvatne@math.ntnu.no}

\begin{document}

\thanks{The authors were supported by Storforsk grant no. 167130 from
  the Norwegian Research Council.}

\maketitle

\begin{abstract}
In this paper we give the derived equivalence classification of
cluster-tilted algebras of type $A_n$. We show that the bounded
derived category of such an algebra depends only on the number of
3-cycles in the quiver of the algebra.
\end{abstract}

\section*{Introduction}

Cluster categories were introduced in \cite{bmrrt} as a framework for
a categorification of Fomin-Zelevinsky cluster algebras
\cite{fz}. In \cite{ccs}, a category was introduced independently for
type $A$, which they showed were equivalent to the cluster
category. For any finite-dimensional hereditary algebra $H$ over a
field $k$, the cluster category $\CC_H$ is the quotient of the bounded derived
category $\DD_H = D^b(\mod H)$ by the functor $F=\tau^{-1}[1]$,
where $\tau$ denotes the AR-translation. $\CC_H$ is canonically
triangulated \cite{keller}, and it has AR-triangles induced by the
AR-triangles in $\DD_H$.

In a cluster category $\CC_H$, tilting objects are defined as objects
which have no self-extensions, and are maximal with respect to this
property. The endomorphism rings of such objects are called
\emph{cluster-tilted algebras} \cite{bmrcta}. These algebras are of
finite representation type if and only if $H$ is the path algebra of a
simply-laced Dynkin quiver.

Cluster-tilted algebras have several interesting properties.
In particular, by \cite{bmrcta} their representation theory can be 
completely understood in terms of the representation theory of the
corresponding 
hereditary algebra $H$. Furthermore,
their relationship to tilted algebras is well understood by \cite{abs,abs2},
see also \cite{ringel}.

Homologically, they are very different from hereditary and tilted algebras,
since they have in general infinite global dimension.
In fact they are 1-Gorenstein and in particular they have finitistic 
dimension 1, by \cite{kr}. Cluster-tilted algebras also play a role in 
the construction of cluster algebras from cluster categories \cite{ck1,ck2},
see also \cite{bmrt}.

The purpose of this paper is to describe when two cluster-tilted
algebras from the cluster category $\CC_H$ have equivalent derived
categories, where $H$ is the path algebra of a quiver whose
underlying graph is $A_n$. We will get an exact description of the
quivers of such algebras, and their relations are given by
\cite{ccs}. The main result is the following.

\begin{Thm}
Two cluster-tilted algebras of type $A_n$ are derived equivalent if
and only if their quivers have the same number of 3-cycles.
\end{Thm}

For this, we show that if we have an almost complete cluster-tilting object
$\ol{T}$ in $\CC_H$ with complements $T_i$ and $T_i^*$ such that the
cluster-tilted algebras given by $\Gamma = \End_{\CC_H}(\ol{T}\amalg
T_i)^{\op}$ and
$\Gamma'= \End_{\CC_H}(\ol{T}\amalg T_i^*)^{\op}$ have quivers with the same
number of 3-cycles, then $\Gamma'$ is in a natural way isomorphic to
the endomorphism ring of a tilting module over $\Gamma$. Then it is
well known that $\Gamma$ and $\Gamma'$ are derived equivalent, see
\cite{happel,cps}.

The outline of the paper is as follows: After some basic notions, we
describe the \emph{mutation class} of $A_n$, that is, the 
quivers of cluster-tilted algebras of $A_n$-type. In Section
\ref{cartdetsec} we give a simple proof of a special case of a result
by Holm \cite{holm}, which is a formula for the determinant of the
Cartan matrices of the cluster-tilted algebras of $A_n$-type. We use
this to distinguish between algebras of this type which are not
derived equivalent. In Section \ref{mainressec} we prove the main
result.

We would like to thank Thorsten Holm and David Smith for interesting
discussions, and Smith for pointing out a missing argument in an
earlier version of the proof of Theorem \ref{mainres}. We would like
to thank Ahmet Seven and an anonymous referee for pointing out missing
references in Sections \ref{mutsec} and \ref{relsec}. There is related
work by G. Murphy \cite{murphy} and Assem, Br\" ustle,
Charbonneau-Jodoin and Plamondon \cite{abcp}.

For notions and basic results about finite dimensional algebras, we
refer the reader to \cite{ass} or \cite{ars}.

\section{Preliminaries}

We will now review some basic notions concerning
cluster-tilted algebras. This theory is developed in
\cite{bmrrt,bmrcta}, and in the Dynkin case there is an independent
approach in \cite{ccs,ccs2}.

Throughout, $H$ will denote the path algebra $k\vec{A}_n$ of a quiver
$\vec{A}_n$ with underlying graph $A_n$. By $\mod H$ we will mean the
category of finitely generated left $H$-modules. Then the AR-quiver of
the derived category $\DD = D^b(\mod H)$ is isomorphic to the stable
translation quiver $\mathbb{Z} A_n$ (see e.g. \cite{happel}). $\DD$
does not depend on the orientation of $\vec{A}_n$.

If $\tau$ is the AR-translation in $\DD$, we consider the
functor $F=\tau^{-1} [1]$ and the orbit category $\CC
=\DD /F$. Then $\CC$ is called the \emph{cluster category} of type
$A_n$. This is a Krull-Schmidt category, and it follows from
\cite{keller} that it has a triangulated structure inherited from $\DD$.

A \emph{(cluster) tilting object} in $\CC$ is an object $T$ with $n$
non-isomorphic indecomposable direct summands such that
$\Ext_{\CC}^{1}(T,T)=0$. An object in $\CC$ with $n-1$ non-isomorphic direct
summands satisfying the same $\Ext$-condition will be called an
\emph{almost complete (cluster) tilting object}. An indecomposable
object $M$ such that $\ol{T}\amalg M$ is a tilting object is said to
be a \emph{complement} of $\ol{T}$.

We will use the following result, which is one of the main results in
\cite{bmrrt}, and which uses the notion of \emph{approximations} from
\cite{as}:

\begin{thm} \label{complement}
An almost complete tilting object $\ol{T}$ in $\CC$ has exactly two
complements $M$ and $M^*$. These are related by unique triangles
$$M\to B\to M^*\to$$
and
$$M^*\to B'\to M\to$$
where the maps $M\to B$ and $M^*\to B'$ are minimal left $\add
\ol{T}$-approximations and the maps $B\to M^*$ and $B'\to M$ are
minimal right $\add \ol{T}$-approximations.
\end{thm}

For a tilting object $T$ in $\CC$, we call the endomorphism ring
$\Gamma_T = \End_{\CC}(T)^{\op}$ a \emph{cluster-tilted
  algebra}. There is a close connection between the module category of
$\Gamma_T$ and $\CC$, from \cite{bmrcta}:

\begin{thm} \label{functor}
With $\Gamma_T$ as above, the functor $G=\Hom_{\CC}(T,-):\CC \to \mod
\Gamma_T$ is full and dense and induces an equivalence
$$\ol{G}:\CC /\add(\tau T) \to \mod \Gamma_T$$
\end{thm}

By \cite{bmrmut}, the cluster-tilted algebras of type $A_n$ are
exactly the algebras given by quivers obtained from $A_n$-quivers by
\emph{mutation}, an operation which will be described in Section
\ref{mutsec}, with certain relations determined by the quiver \cite{bmrfin}.

\section{The mutation class of $A_n$} \label{mutsec}

In this section we will provide an explicit description of the
mutation class of $A_n$-quivers. The ideas underlying our presentation
can be found already in \cite{ccs}, where a geometric interpretation
of mutation of $A_n$-quivers is given. The mutation class is implicit in
\cite{ccs}, see also \cite{seven} for an explicit, but slightly
differently formulated description. The technical Lemma
\ref{iteratemut} will be crucial in the proof of our main theorem in
Section \ref{mainressec}.

\emph{Quiver mutation} was introduced by Fomin and Zelevinsky
\cite{fz} as a generalisation of the sink/source reflections used in
connection with BGP functors \cite{bgp}. Any quiver $Q$ with no
loops and no cycles of length two, can be mutated
at vertex $i$ to a new quiver $Q^*$ by the following rules:
\begin{itemize}
\item The vertex $i$ is removed and replaced by a vertex $i^*$, all
  other vertices are kept.
\item For any arrow $i\to j$ in $Q$ there is an arrow $j\to i^*$ in
  $Q^*$.
\item For any arrow $j\to i$ in $Q$ there is an arrow $i^*\to j$ in
  $Q^*$.
\item If there are $r>0$ arrows $j_1 \to i$, $s>0$ arrows $i\to j_2$
  and $t$ arrows $j_2 \to j_1$ in $Q$, there are $t-rs$ arrows $j_2
  \to j_1$ in $Q^*$. (Here, a negative number of arrows means arrows
  in the opposite direction.)
\item all other arrows are kept
\end{itemize}

Note that if we mutate $Q$ at vertex $i$, and then mutate $Q^*$ at
$i^*$, the resulting quiver is isomorphic to (and will be identified
with) $Q$.
We want to describe the class of quivers which can be obtained by
iterated mutation on a quiver of type $A_n$. Such quivers are said to
be \emph{mutation equivalent} to $A_n$, as iterated mutation produces
an equivalence relation.

The following lemma is a well-known fact:

\begin{lem} \label{nomut}
All orientations of $A_n$ are mutation equivalent.
\end{lem}

From now on, let $\mathcal{Q}_n$ be the class of quivers with $n$
vertices which satisfy the following:
\begin{itemize}
\item all non-trivial cycles are oriented and of length 3
\item a vertex has at most four neighbours
\item if a vertex has four neighbours, then two of its adjacent arrows
  belong to one 3-cycle, and the other two belong to another 3-cycle
\item if a vertex has exactly three neighbours, then two of its
  adjacent arrows belong to a 3-cycle, and the third arrow does not
  belong to any 3-cycle
\end{itemize}
Note that by a \emph{cycle} in the first condition we mean a cycle in
the underlying graph, not passing through the same edge twice. In
particular, this condition excludes multiple arrows. We will show
that $\mathcal{Q}_n$ is the mutation class of $A_n$.

\begin{lem} \label{quivers}
$\mathcal{Q}_n$ is closed under quiver mutation.
\end{lem}

\begin{proof}
Let $Q\in \mathcal{Q}_n$. We will see what happens locally when we mutate.

If we mutate at a vertex $i$ which is a source or a sink, then the
arrows to or from $i$ changes direction, and everything else is left
unchanged. Thus the new quiver $Q^*$ will also satisfy the conditions
in the description of $\mathcal{Q}_n$.

Next we consider the case where $i$ is the source of exactly one arrow
and the target of exactly one arrow:
$$\xymatrix{
j\ar[r] & i\ar[r] & k
}$$
Two cases can occur. Suppose first that there is no arrow from $k$ to
$j$ in $Q$. Then there is an arrow from $j$ to $k$ in $Q^*$:
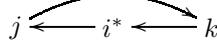
$$\xymatrix{
j\ar @/^1pc/ [rr] & i^*\ar[l] & k\ar[l]
}$$
Thus the numbers of neighbours for $j$ and $k$ increase by 1. It is
impossible that $j$ or $k$ has four neighbours in $Q$, since then the
arrow to or from $i$ would be part of a 3-cycle in $Q$, and $i$ would
have a third neighbour. Thus $j$ and $k$ have $\leq 4$ neighbours in
$Q^*$ as well. There are no other (non-oriented) paths between $j$ and
$k$ in $Q^*$ than the two pictured in the diagram above, so the other
conditions are also satisfied: If $j$ or $k$ has four neighbours in
$Q^*$, then the last two arrows will be part of a 3-cycle in both $Q$
and $Q^*$.

In the other case, there is an arrow $k\to j$ in $Q$. Then this is
removed in passing to $Q^*$. The numbers of neighbours of $j$ and $k$
decrease by 1, and cannot be larger than 3. If, say, $j$ has three
neighbours in $Q^*$, then it must have had four neighbours in $Q$, and
the two arrows not involving $i$ or $k$ are part of a 3-cycle in both
$Q$ and $Q^*$. The arrow $i^*\to j$ is not part of a 3-cycle, since
the only arrow with $i^*$ as target comes from $k$, and there is no
arrow $j\to k$.

We use similar arguments for the other cases, and just point out how
the mutations work. Now let $i$ be a vertex of $Q$ with three
neighbours. Suppose first that the one arrow to or from $i$ which is
not on a 3-cycle has $i$ as the target:
$$
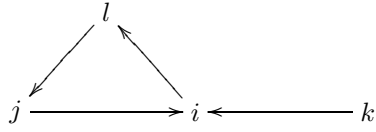
\xymatrix{
&l\ar[ld]&&& \\
j\ar[rr]&&i\ar[ul]&&k\ar[ll]
}$$
Then the mutation will remove the $lj$-arrow and produce a new
triangle $i^*kl$:
$$
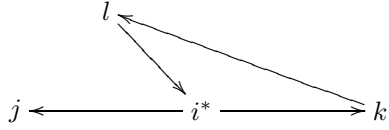
\xymatrix{
&l\ar[rd]&&& \\
j&&i^*\ar[ll]\ar[rr]&&k\ar[ulll]
}$$
Similarly for the case where the third arrow has $i$ as the source.

Finally, let $i$ be a vertex with four neighbours:
$$
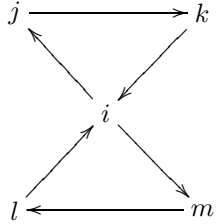
\xymatrix{
j\ar[rr]&&k\ar[ld] \\
&i\ar[ul]\ar[rd]& \\
l\ar[ur]&&m\ar[ll]
}$$
Mutate:
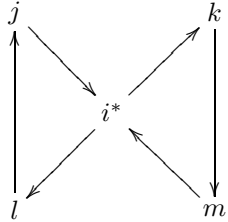
$$\xymatrix{
j\ar[rd]&&k\ar[dd] \\
&i^*\ar[dl]\ar[ru]& \\
l\ar[uu]&&m\ar[ul]
}$$
So for the cases where $i$ has three or four neighbours, we see that
neither in $Q$ nor in $Q^*$ are there other paths between $j,k,l$ and
$m$ than those passing through the diagrams. By similar arguments as
above, $Q^*$ also satisfies the conditions in the description of
$\mathcal{Q}_n$.
\end{proof}

We will need the following lemma for the proof of the main result in
Section \ref{mainressec}.

\begin{lem} \label{iteratemut}
If $Q_1$ and $Q_2$ are quivers in $\mathcal{Q}_n$,
and $Q_1$ and $Q_2$ have the same number of 3-cycles, then $Q_2$ can
be obtained from $Q_1$ by iterated mutation where all the
intermediate quivers also have the same number of 3-cycles.
\end{lem}

\begin{proof}
It is enough to show that all quivers in $\mathcal{Q}_n$ can be mutated without
changing the number of 3-cycles to a quiver looking like this:
$$
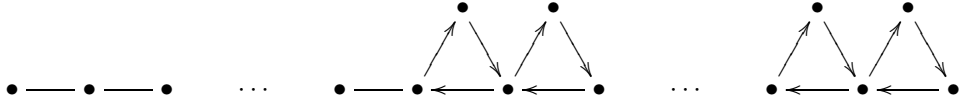
\xymatrix @C0.5pc @R0.6pc{
&&&&&&&&&&&\bullet\ar[ddr]\ar@{<-}[ddl]&&\bullet\ar@{<-}[ddl]\ar[ddr]&&&&&&
\bullet\ar[ddr]\ar@{<-}[ddl]&&\bullet\ar[ddr]\ar@{<-}[ddl]& \\
&&&&&&&&&&&&&&&&&&&&&& \\
\bullet\ar@{-}[rr]&&\bullet\ar@{-}[rr]&&\bullet&&\cdots&&
\bullet\ar@{-}[rr]&&\bullet\ar@{<-}[rr]&&
\bullet&&\bullet\ar[ll]&&\cdots&&\bullet\ar@{<-}[rr]&&
\bullet\ar@{<-}[rr]&&\bullet}$$
In this process we are only allowed to mutate in sinks, sources and
vertices of valency three and four, as these are the mutations which
will not change the number of 3-cycles for quivers in $\mathcal{Q}_n$.

For the purposes of this proof, we introduce a distance function on the
set of 3-cycles in quivers in $\mathcal{Q}_n$. For each pair
$C$, $C'$ of \emph{different} 3-cycles
in $Q$, we define $d_Q(C,C')$ to be the length of the unique minimal
(perhaps non-oriented) path between $C$ and $C'$, i.e. the number
of arrows in this path.

Let $Q$ be a quiver in $\mathcal{Q}_n$, and suppose that the underlying
graph of $Q$ is not $A_n$. We now define a total order on a subset $\SS_Q$ of
the set of 3-cycles of $Q$. This subset is not uniquely defined. $Q$
must contain a 3-cycle
which is only connected to other 3-cycles through (at most) one of its
vertices. Choose one such 3-cycle and call it $C_1$. If there are more
3-cycles, let $C_2$ be
the unique 3-cycle which minimises $d_Q(C_1,-)$. If there are more
3-cycles, let $C_3$ be one of the at most two which minimise
$d_Q(C_2,-)$ among the 3-cycles not equal to $C_1$.

If $C_i$ is defined for some $i\geq 3$, and there exists one or more
3-cycles $C$ such
that $d_Q(C_i,C)<d_Q(C_j,C)$ for $j<i$, let $C_{i+1}$ be one of the at
most two which minimise $d_Q(C_i,-)$ among 3-cycles with this
property. Continue in this way until $C_s$ is defined, but $C_{s+1}$
cannot be defined. Let $\mathcal{S}_Q=\{ C_1,...,C_s\}$ be our totally
ordered set of 3-cycles.

Next, we will see that we have a procedure for moving 3-cycles in the
quiver closer together.
Let $C$ and $C'$ be a pair of neighbouring 3-cycles in $Q$ (i.e. no
edge in the path between them is part of a 3-cycle) such that
$d_Q(C,C')\geq 1$. We want to move $C$ and $C'$ closer together by mutation.
Up to orientation on the arrow from $d$ to $e$, it
looks like the following diagram. The other orientation gives a
similar situation.
$$
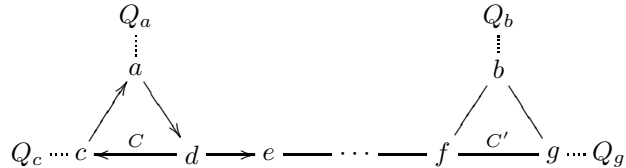
\xymatrix @C0.5pc @R0.5pc{
&&Q_a&&&&&&&&Q_b&& \\
&&a\ar@{.}[u]&&&&&&&&b\ar@{.}[u]&& \\
&&&&&&&&&&&& \\
Q_c&c\ar@{.}[l]\ar[uur]\ar@{<-}[rr]^C&&d\ar@{<-}[uul]\ar[rr]
&&e\ar@{-}[rr]&&\cdots&&f\ar@{-}[ll]\ar@{-}[uur]\ar@{-}[rr]^{C'}
&&g\ar@{-}[uul]\ar@{.}[r]&Q_g
}$$
(In the diagram, the $Q_i$ are subquivers.) Mutating at 
$d$ will produce a quiver $Q^*$ which looks
like this:
$$
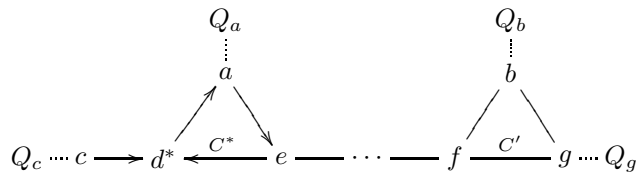
\xymatrix @C0.5pc @R0.5pc{
&&&&Q_a&&&&&&Q_b&& \\
&&&&a\ar@{.}[u]&&&&&&b\ar@{.}[u]&& \\
&&&&&&&&&&&& \\
Q_c&c\ar@{.}[l]\ar[rr]&&d^*\ar[uur]\ar@{<-}[rr]^{C^*}
&&e\ar@{<-}[uul]\ar@{-}[rr]&&\cdots&&f\ar@{-}[ll]\ar@{-}[uur]\ar@{-}[rr]^{C'}
&&g\ar@{-}[uul]\ar@{.}[r]&Q_g
}$$
The only differences between $Q$ and $Q^*$ are that
$d_{Q^*}(C^*,C')=d_Q(C,C')-1$, and there is
after the mutation a path of length 1 between $C^*$ and $Q_c$.

This is the kind of mutation we use for moving 3-cycles closer together.

Suppose that there is a 3-cycle $C$ in $Q$ which is not in our sequence
$\SS_Q$. We will now use the procedure of moving 3-cycles to produce a new
quiver $Q^*$ with a sequence $\mathcal{S}_{Q^*}$ of
3-cycles such that the size of $\SS_{Q^*}$ equals the size of $\SS_Q$
plus one.

The quiver $Q$, with its sequence $\SS_Q$, looks like this:
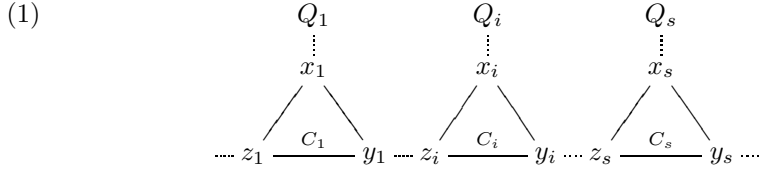
\begin{equation} \label{cyclesequence}
\xymatrix @C=0.5pc @R=0.6pc{
&&Q_1&&&Q_i&&&Q_s&& \\
&&x_1\ar@{.}[u]&&&x_i\ar@{.}[u]&&&x_s\ar@{.}[u]&& \\
&&&&&&&&&& \\
&z_1\ar@{.}[l]\ar@{-}[uur]\ar@{-}[rr]^{C_1}&&y_1\ar@{-}[uul]
&z_i\ar@{.}[l]\ar@{-}[uur]\ar@{-}[rr]^{C_i}&&y_i\ar@{-}[uul]
&z_s\ar@{.}[l]\ar@{-}[uur]\ar@{-}[rr]^{C_s}&&y_s\ar@{-}[uul]\ar@{.}[r]&
}
\end{equation}
where the $Q_i$ are subquivers, and $C$ is in $Q_i$ for some
$i=2,3,...,s-1$. (This follows from the definition of $C_1$ and $s$.)
$C$ may be moved towards $x_i$ 
using the procedure above. So we may perform this procedure until $C$
and $C_i$ share the vertex $x_i$, and $x_i$ has four
neighbours. We then mutate at the vertex $x_i$:
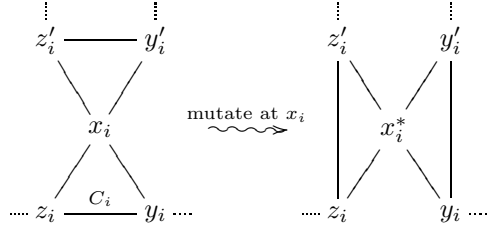
$$
\xymatrix @C=0.5pc @R=0.5pc{
&&&&&&&&&&& \\
&z_i'\ar@{.}[u]\ar@{-}[rr]&&y_i'\ar@{.}[u]&&&&&
z_i'\ar@{.}[u]&&y_i'\ar@{.}[u]& \\
&&&&&&&&&&& \\
&&x_i\ar@{-}[uur]\ar@{-}[uul]\ar@{-}[ddr]\ar@{-}[ddl]&&\ar@{~>}[rrr]^{\textrm{mutate
  at }x_i}&&&&&
x_i^*\ar@{-}[uur]\ar@{-}[uul]\ar@{-}[ddl]\ar@{-}[ddr]&& \\
&&&&&&&&&&& \\
&z_i\ar@{.}[l]\ar@{-}[rr]^{C_i}&&y_i\ar@{.}[r]&&&&&
z_i\ar@{-}[uuuu]\ar@{.}[l]&&y_i\ar@{-}[uuuu]\ar@{.}[r]&
}
$$
Call the resulting quiver $Q^*$. After a suitable labelling, we now have a
sequence $C_1^*,...,C_{s+1}^*$ of 3-cycles in the quiver $Q^*$, where
$C_j^*=C_j$ for $j<i$ and $C_j^*=C_{j-1}$ for $j>i+1$. This may
serve as a sequence $\SS_{Q^*}$.

Enlarging our totally ordered set like this the necessary number of
times will give a quiver where
\emph{all} the 3-cycles are in a sequence $C_1,...,C_s$ as in 
diagram \eqref{cyclesequence} for some $s$, and the subquivers
$Q_1,...,Q_s$ are just (non-directed) paths.

If $y_s$ in diagram \eqref{cyclesequence} has valency 3, we now move
$C_s$ to the right by mutating at $y_s$ and continuing in the
same way. When we reach a diagram as in \eqref{cyclesequence} above
where $y_s$ has only two
neighbours ($x_s$ and $z_s$), we shrink $Q_s$ by mutating at $x_s$ and
continuing until the new $x_s$ has only $y_s$ and $z_s$ as
neighbours. By suitably orienting $Q_s$ beforehand as in Lemma
\ref{nomut}, we can do this in
such a way that $y_s$ still only has two neighbours, and $C_s$ is
connected to the rest of the quiver only through $z_s$. Successively
doing this to $C_{s-1},...,C_1$ will give a quiver consisting of a
sequence of 3-cycles with $d_Q(C_i,C_{i+1})=0$ for neighbouring $C_i$
and $C_{i+1}$, and possibly with some non-directed path connected to it:
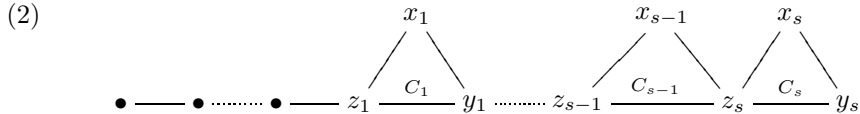
\begin{equation} \label{standardquiver}
\xymatrix @C=0.5pc @R=0.6pc{
&&&&&&&x_1\ar@{-}[ddl]\ar@{-}[ddr]&&&&x_{s-1}\ar@{-}[ddl]\ar@{-}[ddr]
&&x_s\ar@{-}[ddl]\ar@{-}[ddr]& \\
&&&&&&&&&&&&&& \\
\bullet\ar@{-}[rr]&&\bullet\ar@{.}[rr]&&\bullet\ar@{-}[rr]
&&z_1\ar@{-}[rr]^{C_1}&&y_1\ar@{.}[rr]&&z_{s-1}\ar@{-}[rr]^{C_{s-1}}
&&z_s\ar@{-}[rr]^{C_s}&&y_s
}
\end{equation}

The orientation of $C_s$ does not matter, since we can just flip it in
the diagram. If $i$ is the biggest number $<s$ such that $C_i$ is not
oriented in the clockwise direction, we mutate at $y_i=z_{i+1}$ and
get a similar diagram where the new $C_{i+1}$ is oriented in the
anticlockwise direction, and the new $C_i$ is oriented
clockwise. Doing this the necessary number of times, we get the quiver
we want.
\end{proof}

It should be remarked that the following proposition follows from
Lemma \ref{quivers} and the fact that quivers in $\mathcal{Q}_n$ are
2-finite \cite{fz2}, see also \cite{seven}. However, we can now give
an independent argument:

\begin{prop}
A quiver $Q$ is mutation equivalent to $A_n$ if and only if $Q\in
\mathcal{Q}_n$.
\end{prop}

\begin{proof}
Obviously, all orientations of $A_n$ are in $\mathcal{Q}_n$.

It follows from the proof of Lemma \ref{iteratemut} that all members
of $\mathcal{Q}_n$ can be reached by iterated mutation on an
$A_n$-quiver, since mutating the quiver in \eqref{standardquiver} in
all the $x_i$ will give a quiver with underlying graph $A_n$, and we
can reverse the procedure in the proof to come to any $Q\in
\mathcal{Q}_n$.
\end{proof}

\section{Relations} \label{relsec}

In this section we give the relations on the quivers of cluster-tilted
algebras of type $A_n$, which are given in \cite{ccs} and have
been generalised in \cite{ccs2} and \cite{bmrfin}. This gives the complete
description of this class of algebras, and we use it to establish that
these algebras are gentle.

\begin{prop}
The cluster-tilted algebras of type $A_n$ are exactly the algebras
$kQ/I$ where $Q$ is a quiver in $\mathcal{Q}_n$ and $I$ is the ideal
generated by the directed paths of length 2 which are part of a 3-cycle.
\end{prop}

Given such a quiver $Q$, we will sometimes denote the
corresponding cluster-tilted algebra $kQ/I$ by $\Gamma_Q$.

If $Q$ is a finite quiver and $I$ is an ideal in the path algebra $kQ$,
then $kQ/I$ is \emph{special biserial} \cite{skw} if it satisfies
\begin{itemize}
\item For every vertex $p$ in $Q$, there are at most two arrows
  starting in $p$ and at most two arrows ending in $p$.
\item For every arrow $\beta$ in $Q$, there is at most one arrow
  $\alpha_1$ in $Q$ with $\beta \alpha_1\notin I$ and at most one arrow
  $\gamma_1$ in $Q$ with $\gamma_1 \beta\notin I$.
\end{itemize}

A special biserial algebra $kQ/I$ is \emph{gentle} \cite{assk} if it
also satisfies
\begin{itemize}
\item $I$ is generated by paths of length 2.
\item For every arrow $\beta$ in $Q$ there is at most one arrow
  $\alpha_2$ such that $\beta \alpha_2$ is a path and $\beta
  \alpha_2\in I$, and at most one arrow $\gamma_2$ such that $\gamma_2
  \beta$ is a path and $\gamma_2 \beta\in I$.
\end{itemize}

\begin{cor}
Cluster-tilted algebras of type $A_n$ are gentle.
\end{cor}

\section{Cartan determinants} \label{cartdetsec}

The \emph{Cartan matrix} $(C_{ij})$ of a finite-dimensional
$k$-algebra $\Lambda$ is by definition the matrix with $ij$'th
entry $C_{ij}=\dim_k e_i\Lambda e_j$, that is, the columns are the
dimension vectors of the indecomposable projectives. The determinant
of the Cartan matrix is invariant under derived equivalence. (See
\cite{bocianskow} for a proof.)

Since cluster-tilted algebras of type $A_n$ are gentle,
the following result is a special case of a result by Holm
\cite{holm}. We include the proof, which is a lot simpler than in the
general case.

\begin{prop} \label{cartdet}
If $\Gamma$ is a cluster-tilted algebra of type $A_n$, then the
determinant of the Cartan matrix of $\Gamma$ is $2^{t}$, where $t$ is
the number of 3-cycles in the quiver of $\Gamma$.
\end{prop}

\begin{proof}
Let $Q$ be any quiver with relations $R$, and let $\Lambda_{Q,R}$ be
the corresponding algebra. We will study the behaviour of the Cartan
determinant of the algebra under two types of enlargements of $Q$ and
$R$. All quivers in the mutation class of $A_n$ can be built by
successive enlargements of this kind, and this will yield the result.

Let $k$ denote the number of vertices in $Q$.

The first type of enlargement goes as follows.
Construct $Q'$ from $Q$ by adding one new vertex labelled $k+1$ and
one arrow $\alpha$ between $Q$ and k+1. We assume that $k+1$
is the target of $\alpha$. (The proof is similar in the other
case.) Without loss of generality we can also assume that the source of
$\alpha$ is the vertex labelled $k$. Let the relations $R'$ on $Q'$ be
$R$, i.e. $Q'$ just inherits the relations from $Q$. Then the
Cartan matrix $C(\Lambda_{Q',R'})$ is
$$\left(
\begin{array}{cccc}
*&\cdots&*&0 \\
\vdots&&\vdots&\vdots \\
*&\cdots&*&0 \\
*&\cdots&*&1
\end{array}\right)
$$
where the Cartan matrix $C(\Lambda_{Q,R})$ sits in the top left
corner. We see that the determinant of $C(\Lambda_{Q',R'})$ equals the
determinant of $C(\Lambda_{Q,R})$, so this construction does not
change the Cartan determinant.

We now turn to the second type of enlargement. Construct a quiver
$Q''$ from $Q$ by adding two vertices $k+1$ and
$k+2$, and three arrows $\alpha$, $\beta$, $\gamma$ such that $\gamma
\beta \alpha$ is a 3-cycle running through $k+1$ and $k+2$. Again,
we may assume that the third vertex in this 3-cycle is labelled
$k$. Let the relations on $Q''$ be given by $R''=R\cup \{ \beta \alpha,
\gamma \beta, \alpha \gamma\}$.

Since the minimal relations involving $\alpha$, $\beta$ and $\gamma$
do not involve the relations in $R$, we have, for $i\leq k$, 
\begin{eqnarray*}
\dim_k e_i(\Lambda_{Q'',R''})e_{k+1} & = & \dim_k e_i(\Lambda_{Q,R})e_k \\
\dim_k e_{k+2}(\Lambda_{Q'',R''})e_i & = & \dim_k e_k(\Lambda_{Q,R})e_i
\end{eqnarray*}
This gives the following $(k+2)\times (k+2)$ Cartan matrix of
$C(\Lambda_{Q'',R''})$.
$$\left(
\begin{array}{cccccc}
*&\cdots&*&|&0&| \\
\vdots&&\vdots&|&\vdots&| \\
*&\cdots&*&|&0&| \\
-&-&-&1&0&1 \\
-&-&-&1&1&0 \\
0&\cdots&0&0&1&1
\end{array}\right)$$
Here, rows $k$ and $k+1$ are equal except for the two rightmost
entries, and similarly for columns $k$ and $k+2$. Again, we find
$C(\Lambda_{Q,R})$ in the top left corner. Expanding the determinant
along the bottom row, we get
$$|C(\Lambda_{Q'',R''})|=M_{k+2,k+2}-M_{k+2,k+1}$$
where $M_{ij}$ denotes the $ij$'th minor. We see that
$M_{k+2,k+2}=|C(\Lambda_{Q,R})|$, while
\begin{equation*}
M_{k+2,k+1}=
\left|
\begin{array}{ccccc}
*&\cdots&*&|&| \\
\vdots&&\vdots&|&| \\
*&\cdots&*&|&| \\
-&-&-&1&1 \\
-&-&-&1&0
\end{array}
\right|
\end{equation*}
where rows $k$ and $k+1$ are equal except for the last entry, and
similarly for columns $k$ and $k+1$. The top left part is still
$C(\Lambda_{Q,R})$. Upon subtracting row $k+1$ from
row $k$ and then expanding the determinant along row $k$, we find that
$$M_{k+2,k+1}=-|C(\Lambda_{Q,R})|$$
and thus
$$|C(\Lambda_{Q'',R''})|=2\cdot |C(\Lambda_{Q,R})|$$

The statement in the proposition now follows from the following
observation: Given any cluster-tilted algebra $\Gamma$ of type $A_n$,
we may build the quiver of $\Gamma$ with the appropriate relations by
starting with $A_1$ and performing the above types of enlargements
sufficiently many times. The determinant of the Cartan matrix is
multiplied by two for each 3-cycle added.
\end{proof}

\section{Derived equivalence} \label{mainressec}

We now prove the main result, namely that the cluster-tilted algebras
of type $A_n$ which have the same Cartan determinant are also derived
equivalent.

\begin{thm} \label{mainres}
Two cluster-tilted algebras of type $A_n$ are derived equivalent if
and only if their quivers have the same number of 3-cycles.
\end{thm}

\begin{proof}
Let $\Gamma$ and $\Gamma'$ be two cluster-tilted algebras of type
$A_n$.

If the quivers of $\Gamma$ and $\Gamma'$ do not have the same number
of 3-cycles, then by Proposition \ref{cartdet} the determinants of
their Cartan matrices are not equal, and thus they are not derived
equivalent.

By the results in Section \ref{mutsec}, and in particular Lemma
\ref{iteratemut}, it is enough to show that if $\Gamma$ and
$\Gamma '$ are two cluster-tilted algebras of type $A_n$, and the
quiver of one of them can be obtained by mutating the quiver of the other in
one vertex without changing the number of 3-cycles, then $\Gamma$ and
$\Gamma '$ are derived equivalent. The strategy is to show that
replacing a certain direct summand of the tilting object $T$ in $\CC$
with $\Gamma = \End_{\CC}(T)$, we find a tilting $\Gamma$-module whose
endomorphism ring is isomorphic to $\Gamma'$.

If the vertex $i$ for the necessary mutation is a source or a sink,
the mutation just corresponds to APR-tilting, so this is a well-known
case \cite{apr}.

We consider the case where $i$ is a vertex of the quiver of $\Gamma$
with two arrows ending there. There might be one or two arrows with
$i$ as the initial vertex. (Having proved the result for this case, we
need not do it for the case where $i$ is a vertex with two arrows out
and one arrow in, since this is just the reverse operation of what we
have done.)

$$
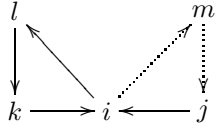
\xymatrix{l \ar[d] & & m \ar@{.>}[d] \\
k \ar[r] & i \ar[ul] \ar@{.>}[ur] & j \ar[l]}$$

Let $T$ be the tilting object in the cluster category $\CC$ which
gives rise to $\Gamma =\End_{\CC}(T)^{\op}$, and suppose $T_i$ is the
indecomposable summand which through the functor
$\Hom_{\CC}(T,-):\CC\to \mod \Gamma$ corresponds to the vertex where
we must mutate to get the quiver of $\Gamma'$ from the quiver of $\Gamma$.
Denote by $T_i^*$ the unique second indecomposable object which
completes $\ol{T}$ to a tilting object in $\CC$. Then
$$\Gamma'=\End_{\CC}(\ol{T}\amalg T_i^*)^{\op}.$$

By Theorem \ref{functor}, the functor $G=\Hom_{\CC}(T,-): \CC\to
\mod \Gamma$ is full and dense and induces an equivalence
$$\ol{G}:\CC /\add (\tau T)\to \mod \Gamma.$$
Let $T_1, ...,T_n$ be the
indecomposable summands of $T$. The images of these objects under
$\Hom_{\CC}(T,-)$ are the
indecomposable projective $\Gamma$-modules:
$$P_i=G(T_i)=\Hom_{\CC}(T,T_i)$$
Denote by $P_i^*$ the image of $T_i^*$:
$$P_i^*=\Hom_{\CC}(T,T_i^*)$$

$T_i^*$ is found by completing a minimal left $\add
\ol{T}$-approximation $T_i\to B$ into a triangle in $\CC$ (cf. Theorem
\ref{complement}):
\begin{equation} \label{triangle}
T_i\to B\to T_i^*\to
\end{equation}
We see that $B$ in this triangle is $T_j\amalg T_k$, where $j$ and $k$
are the labels on the vertices which have arrows to $i$ in the quiver
of $\Gamma$. Recall that the AR-quiver of $\CC$ is a M\" obius band
with $\mod H \vee H[1]$ as a fundamental domain for the functor $F=\tau^{-1}[1]$,
which takes the $\mathbb{Z} A_n$ to the M\" obius band
\cite{bmrrt}. The AR-quiver is drawn in the following diagram, where
the dotted lines indicate a choice of fundamental domain for $F$.
{\tiny
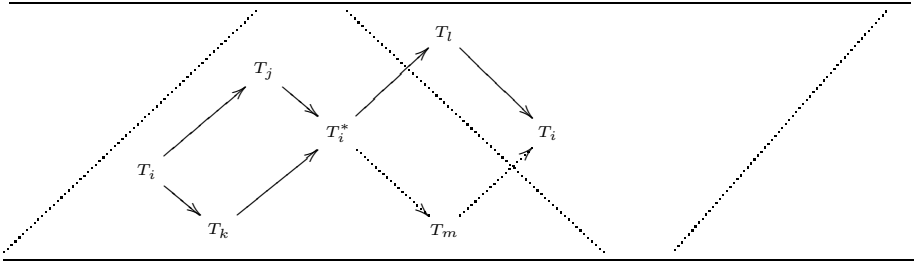
$$\xymatrix @R=0.2pc @C=0.2pc{
\;\ar@{-}[rrrrrrrrrrrrrrrrrrrrrrrrrrrrr]
\;&\;&\;&\;&\;&\;&\;&\;&\;&\;
&\;\ar@{.}[ddddddddrrrrrrrr]&\;&\;&\;&\;&\;&\;
&\;&\;&\;&\;&\;&\;&\;&\;
&\;&\;&\;&\;&\; \\
&&&&&&&&&&&&&{}\save[]T_l\ar[dddrrr]\restore&&&&&&&&&&&&&&&&& \\
&&&&&&&&{}\save[]T_j\ar[ddrr]\restore&&&&&&&&&&&&&&&&&&&&&& \\
&&&&&&&&&&&&&&&&&&&&&&&&&&&&&& \\
&&&&&&&&&&{}\save[]T_i^*\ar[uuurrr]\ar@{.>}[dddrrr]\restore
&&&&&&{}\save[]T_i\restore&&&&&&&&&&&&&& \\
&&&&&{}\save[]T_i\ar[uuurrr]\ar[ddrr]\restore&&&&&&&&&&&&&&&&&&&&&&&&& \\
&&&&&&&&&&&&&&&&&&&&&&&&&&&&&& \\
&&&&&&&
{}\save[]T_k\ar[uuurrr]\restore&&&&&&{}\save[]T_m\ar@{.>}[uuurrr]\restore
&&&&&&&&&&&&&&&&& \\
\ar@{-}[rrrrrrrrrrrrrrrrrrrrrrrrrrrrr]\ar@{.}[uuuuuuuurrrrrrrr]
&&&&&&&&&&&&&&&&&&&& \ar@{.}[uuuuuuuurrrrrrrr]&&&&&&&&&&
}
$$
}

We must show that $P_i^*\neq 0$, which is the same as showing that
$T_i^*$ is not in $\add(\tau T)$. Since $T_i^*$ is indecomposable, this
would mean that $T_i^*=\tau T_q$ for some $q=1,2,...,n$. But by Serre
duality this would lead to
$$\Ext_{\CC}^{1}(T_q,T_j)\simeq D\Hom_{\CC}(T_j,\tau
T_q)=D\Hom_{\CC}(T_j,T_i^*)\neq 0,$$
which is absurd, since $T$ is a tilting object in $\CC$. Thus we
conclude that $P_i^*\neq 0$.

We will now show that no non-zero endomorphisms of $\ol{T}\amalg
T_i^*$ factor through $\add(\tau T)$. Then $G$ and the induced functor
$\ol{G}$ will give the isomorphism
\begin{equation} \label{endoiso}
\Gamma'=\End_{\CC}(\ol{T}\amalg T_i^*)\simeq
\End_{\Gamma}(\ol{P}\amalg P_i^*),
\end{equation}
$$f\mapsto \Hom_{\CC}(T,f).$$
It is only necessary to consider maps involving $T_i^*$. Also, we see
that if a non-zero map factors $T_i^*\to \tau T_q \to T_s$, then
$q=i$, for this implies the existence of an extension between $T_i^*$
and $T_q$. We find the
maps in $\End_{\CC}(\ol{T}\amalg T_i^*)$ from the AR-quiver of $\CC$,
and we recall from \cite{bmrrt} that the $k$-dimension of
$\Hom_{\CC}(X,Y)$ is at most 1 when $X$ and $Y$ are summands of a
tilting object.

We consider first the situation where $i$ is the initial vertex of two
arrows:
{\tiny
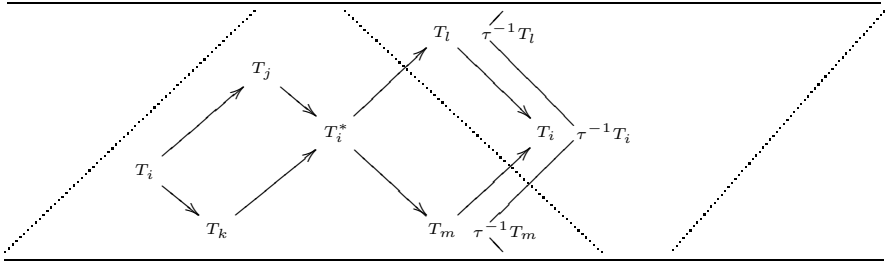
$$\xymatrix @R=0.2pc @C=0.2pc{
\;\ar@{-}[rrrrrrrrrrrrrrrrrrrrrrrrrrrr]
&\;&\;&\;&\;&\;&\;&\;&\;
&\;&\;\ar@{.}[ddddddddrrrrrrrr]&\;&\;&\;&\;&\;
&\;&\;&\;&\;&\;&\;&\;&\;
&\;&\;&\;&\;&\; \\
&&&&&&&&&&&&&{}\save[]T_l\ar[dddrrr]\restore&{}\save[]+<.3cm,0cm>*\txt<8pc>
{$\tau^{-1}T_l$}\restore\ar@{-}[dddrrr]\ar@{-}[ur]
&&&&&&&&&&&&&& \\
&&&&&&&&{}\save[]T_j\ar[ddrr]\restore&&&&&&&&&&&&&&&&&&&&& \\
&&&&&&&&&&&&&&&&&&&&&&&&&&&&& \\
&&&&&&&&&&{}\save[]T_i^*\ar[uuurrr]\ar[dddrrr]\restore
&&&&&&T_i&{}\save[]+<.3cm,0cm>*\txt<8pc>{$\tau^{-1}T_i$}\restore&&&&&&&&&&& \\
&&&&&{}\save[]T_i\ar[uuurrr]\ar[ddrr]\restore&&&&&&&&&&&&&&&&&&&&&&&& \\
&&&&&&&&&&&&&&&&&&&&&&&&&&&&& \\
&&&&&&&{}\save[]T_k\ar[uuurrr]\restore&&&&&&{}\save[]T_m\ar[uuurrr]\restore
&{}\save[]+<.3cm,0cm>*\txt<8pc>{$\tau^{-1}T_m$}\restore\ar@{-}[uuurrr]\ar@{-}[dr]
&&&&&&&&&&&&&&& \\
\ar@{-}[rrrrrrrrrrrrrrrrrrrrrrrrrrrr]\ar@{.}[uuuuuuuurrrrrrrr]
&&&&&&&&&&&&&&&&&&&& \ar@{.}[uuuuuuuurrrrrrrr]&&&&&&&&&
}
$$
}
Suppose there is a non-zero map $\phi : T_i^*\to T_s$ for some
$s$. Then $\phi$ factors through either $T_l$ or $T_m$, since the map
$T_i^*\to T_l\amalg T_m$ is a minimal left
$\add(\ol{T})$-approximation. Also, $\phi$ cannot factor through
$\tau^{-1}T_l$ or $\tau^{-1}T_m$ since $\ol{T}$ is exceptional.
Therefore $T_s$ must be on one of the rays starting in
either $T_l$ or $T_m$. Since there are arrows in the quiver of
$\Gamma$ from $i$ to $l$ and $m$, $T_s$ cannot be on the rays from
$T_l$ and $T_m$ to $T_i$ in the AR-quiver, so it must be on one of the
rays starting in $T_i^*$, after $T_l$ or $T_m$.

So the only way $\phi$ can factor through $\tau T_i$ is if one of
the maps $T_l\to T_i$ and $T_m\to T_i$ is irreducible. But this is
impossible, since this would imply that there is a non-zero map
$\tau^{-1}T_l\to T_k$, i.e. an extension between $T_l$ and $T_k$, or
similarly for $T_m$ and $T_j$.

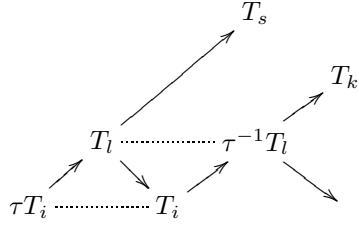
$$\xymatrix @R=.8pc @C=.8pc{
&&&{}\save[]T_s\restore& \\
&&&&{}\save[]T_k\restore \\
&{}\save[]T_l \ar[uurr] \ar[dr]
\ar@{.}[rr]\restore&&{}\save[]\tau^{-1} T_l\ar@{->}[ur] \ar@{->}[dr]\restore& \\
{}\save[]\tau T_i \ar@{.}[rr] \ar@{->}[ur]\restore&&{}\save[]T_i
\ar@{->}[ur]\restore&&
}$$

Now consider the case where $i$ is the initial vertex of only one arrow.
{\tiny
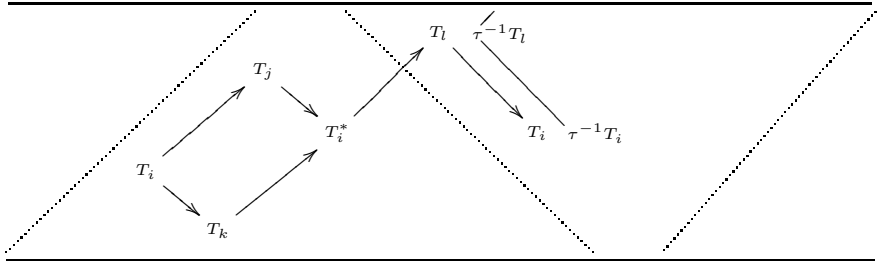
$$\xymatrix @R=0.2pc @C=0.2pc{
\;\ar@{-}[rrrrrrrrrrrrrrrrrrrrrrrrrrrr]
&\;&\;&\;&\;&\;&\;&\;&\;&\;
&\;\ar@{.}[ddddddddrrrrrrrr]&\;&\;&\;&\;&\;
&\;&\;&\;&\;&\;&\;&\;&\;
&\;&\;&\;&\;&\; \\
&&&&&&&&&&&&&{}\save[]T_l\ar[dddrrr]\restore&{}\save[]+<.3cm,0cm>*\txt<8pc>
{$\tau^{-1}T_l$}\restore\ar@{-}[dddrrr]\ar@{-}[ur]
&&&&&&&&&&&&&& \\
&&&&&&&&{}\save[]T_j\ar[ddrr]\restore&&&&&&&&&&&&&&&&&&&&& \\
&&&&&&&&&&&&&&&&&&&&&&&&&&&&& \\
&&&&&&&&&&{}\save[]T_i^*\ar[uuurrr]\restore
&&&&&&T_i&{}\save[]+<.3cm,0cm>*\txt<8pc>{$\tau^{-1}T_i$}\restore&&&&&&&&&&& \\
&&&&&{}\save[]T_i\ar[uuurrr]\ar[ddrr]\restore&&&&&&&&&&&&&&&&&&&&&&&& \\
&&&&&&&&&&&&&&&&&&&&&&&&&&&&& \\
&&&&&&&{}\save[]T_k\ar[uuurrr]\restore&&&&&&&
&&&&&&&&&&&&&&& \\
\ar@{-}[rrrrrrrrrrrrrrrrrrrrrrrrrrrr]\ar@{.}[uuuuuuuurrrrrrrr]
&&&&&&&&&&&&&&&&&&&& \ar@{.}[uuuuuuuurrrrrrrr]&&&&&&&&&
}
$$
}
Again, a non-zero map $\phi : T_i^*\to T_s$ must factor through $T_l$,
since $T_i^*\to T_l$ is a minimal left
$\add(\ol{T})$-approximation. For the same reasons as above, it cannot
factor through $\tau^{-1}T_l$, so $T_s$ must be on one of the rays
starting in $T_l$, and it can not be between $T_l$ and $T_i$. The
possibility that it is on the same ray as $T_i$ is then ruled out by
the fact that the map $T_i^*\to T_l\to T_i$ is a composition of two maps in a
triangle and therefore zero. This forces us to the situation described
above with an irreducible map $T_l\to T_i$, which is again
impossible. Summarising, we have that no non-zero morphisms $T_i^*\to
T_s$ can factor through a $\tau T_q$.

Next note that a non-zero morphism $T_s \to T_i^*$ factoring through a
$\tau T_q$ would provide an extension between $T_q$ and $T_s$, which
is not possible. Obviously, there cannot be any endomorphisms of
$T_i^*$ factoring through other indecomposable objects. Hence, the
isomorphism \eqref{endoiso} is established.

Our goal is now to show that $\ol{P}\amalg P_i^*$ is a tilting module
over $\Gamma$. Then the isomorphism \eqref{endoiso} will imply that
the derived categories $D^b(\mod \Gamma )$ and $D^b(\mod \Gamma')$ are
equivalent, by Happel's theorem \cite{happel,cps}.

We now claim that in $\mod \Gamma$ the triangle \eqref{triangle} gives
us a projective resolution of $P_i^*$:
\begin{equation} \label{projres}
0\to P_i\to P_j\amalg P_k\to P_i^*\to 0
\end{equation}
So the $\Gamma$-module $P_i^*$ has projective dimension 1. Indeed, the
triangle \eqref{triangle} provides the long exact sequence
\begin{equation} \label{projresII}
\cdots \to \Hom_{\CC}(T,\tau^{-1}T_j\amalg\tau^{-1}T_k)\to
\Hom_{\CC}(T,\tau^{-1}T_i^*)\to
\end{equation}
$$
\Hom_{\CC}(T,T_i)\to\Hom_{\CC}(T,T_j\amalg T_k)\to \cdots
$$
and the map $P_i\to P_j\amalg P_k$ in \eqref{projres} is mono if and
only if the first map in \eqref{projresII} is epi. To see that this is the
case, we consider maps $\psi: \tau T_s\to T_i^*$, $1\leq s\leq n$, and
show that they factor through $T_j\amalg T_k$. This is easily seen
when $s=i$. If $\psi$ factors
through $T_s$, it will also factor through $T_j\amalg T_k\to T_i^*$,
since this is an $\add \ol{T}$-approximation. If $\psi$ does not
factor through $T_s$, $\tau T_s$ must be on one of the rays pointing
to $T_i^*$ in the AR-quiver. Without loss of generality, assume
that this is the ray which $T_k$ lies on. $\tau T_s$ cannot be between
$T_k$ and $T_i^*$, since this would imply an extension between $T_k$ and $T_s$:
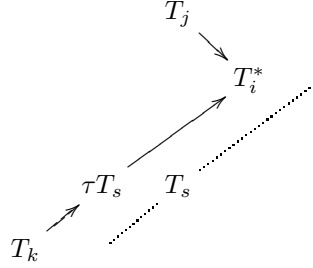
$$\xymatrix @R=.8pc @C=.8pc{
&&&{}\save[]T_j\ar[dr]\restore&& \\
&&&&{}\save[]T_i^*\restore&\quad \\
&&&&& \\
&&{}\save[]\tau T_s\ar[uurr]\restore&{}\save[]T_s\ar@{.}[uurr]\restore&& \\
&{}\save[]T_k\ar[ur]\restore&{}\save[]\ar@{.}[ur]\restore&&&
}$$
In the case where $T_k$ is between $\tau T_s$ and $T_i^*$, $\psi$
factors through $T_k$, which was what we wanted to show. Thus the
first map in \eqref{projresII} is epi, and \eqref{projres} is the
projective resolution of $P_i^*$.

We need to check that $\Ext_{\Gamma}^{a}(P_i^*,P_l)=0$ for $a=1,2,3,...$
and all indecomposable projectives $P_l\neq P_i$. Since $P_i^*$ has
projective dimension 1, the $a\geq 2$ case is trivial. Passing
to the derived category $\DD_{\Gamma} =D^{b}(\mod \Gamma)$, we
identify $P_i^*$ with the deleted projective resolution
$$\cdots \to 0\to P_i\to P_j\amalg P_k\to 0\to \cdots$$
where $P_j\amalg P_k$ sits in degree 0. This derived category is
equivalent to the homotopy category $K^{-,b}(\textrm{proj }\Gamma)$ of
upper bounded complexes of projective modules with non-zero homology only in a
finite number of positions. Since
$$\Ext_{\Gamma}^{1}(X,Y)\simeq \Hom_{\DD_{\Gamma}}(X,Y[1])$$
for objects $X$ and $Y$ of $\DD_{\Gamma}$, we need to show that there are no
non-zero (up to homotopy) morphisms of complexes
$$\xymatrix{\cdots \ar[r] & 0 \ar[r]\ar[d] & P_i \ar[r]\ar[d] & P_j\amalg P_k
  \ar[r]\ar[d] & 0 \ar[r]\ar[d] & \cdots \\
\cdots \ar[r] & 0 \ar[r] & P_q \ar[r] & 0 \ar[r] & 0 \ar[r] & \cdots
}$$
But all morphisms $P_i\to P_q$ where $P_q$ is an indecomposable
projective factor through $P_i \to P_j\amalg P_k$, since this map is a
minimal left $\add \ol{P}$-approximation, so all morphisms of
complexes as above are null-homotopic. Thus all the $\Ext$-groups
vanish.

Since the number of indecomposables in $\ol{P}\amalg P_i^*$ equals $n$, the
number of simples over $\Gamma$, and $\ol{P}\amalg P_i^*$ has
projective dimension 1, it is a tilting module, and we are done.
\end{proof}


\begin{thebibliography}{99}

\bibitem[ABCP]{abcp} Assem I., Br\" ustle T., Charbonneau-Jodoin G.,
  Plamondon P.-G. \emph{Cluster-tilted gentle algebras}, in preparation

\bibitem[ABS1]{abs} Assem I., Br\"ustle T., Schiffler R. 
\emph{Cluster-tilted algebras as trivial extensions}, To appear in
J. London Math. Soc., preprint arxiv:math.RT/0601537

\bibitem[ABS2]{abs2} Assem I., Br\" ustle T., Schiffler
  R. \emph{Cluster-tilted algebras and slices}, preprint v.2
  arxiv:math.RT/0707.0038 (2007)

\bibitem[AsSk]{assk} Assem I., Skowro\' nski A. \emph{Iterated tilted
    algebras of type $\tilde{A}_n$}, Math. Z. 195, 269-290 (1987)

\bibitem[ASS]{ass} Assem I., Simson D., Skowro\' nski
  A. \emph{Elements of the Representation Theory of Associative
    Algebras Vol. 1: Techniques of representation theory}, London
  Mathematical Society Student Texts 65, Cambridge University Press
  (2006)

\bibitem[APR]{apr} Auslander M., Platzeck M. I., Reiten
  I. \emph{Coxeter functors without diagrams},
  Trans. Amer. Math. Soc. 250, 1-46 (1979)

\bibitem[ARS]{ars} Auslander M., Reiten I., Smal\o\
  S. \emph{Representation Theory of Artin Algebras}, Cambridge Studies
  in Advanced Mathematics 36, Cambridge University Press (1997)

\bibitem[AS]{as} Auslander M., Smal\o\ S. \emph{Preprojective
    modules over Artin algebras}, J. Algebra 66, no. 1, 61-122 (1980)

\bibitem[BGP]{bgp} Bernstein I. N., Gelfand I. M., Ponomarev
  V. A. \emph{Coxeter functors, and Gabriel's theorem}, Uspehi
  Mat Nauk 28, no. 2, 19-33 (1973)

\bibitem[BoSk]{bocianskow} Bocian R., Skowro\' nski \emph{Weakly
    symmetric algebras of Euclidean type}, J. Reine Angew. Math. 580
  157-199 (2005)

\bibitem[BMR1]{bmrcta} Buan A., Marsh R., Reiten I. \emph{Cluster-tilted
    algebras}, Trans. Amer. Math. Soc. 359, no. 1, 323-332 (2007)

\bibitem[BMR2]{bmrmut} Buan A., Marsh R., Reiten I. \emph{Cluster
    mutation via quiver representations}, preprint v.2 math.RT/0412077
  (2005)

\bibitem[BMR3]{bmrfin} Buan A., Marsh R., Reiten
  I. \emph{Cluster-tilted algebras of finite representation type},
  J. Algebra 306, no. 2, 412-431 (2006)

\bibitem[BMRRT]{bmrrt} Buan A., Marsh R., Reineke M., Reiten I.,
  Todorov G. \emph{Tilting theory and cluster combinatorics}, Adv.
  Math. 204, 572-618 (2006)

\bibitem[BMRT]{bmrt} Buan A., Marsh R., Reiten I., Todorov G. 
\emph{Clusters and seeds for acyclic cluster algebras} with an appendix by
Buan A., Caldero P., Keller B., Marsh R., Reiten I., Todorov G., 
Proc. Amer. Math. Soc. 135, no. 10, 3049-3060 (2007)

\bibitem[CCS1]{ccs} Caldero P., Chapoton F., Schiffler R. \emph{Quivers
    with relations arising from clusters ($A_n$ case)},
  Trans. Amer. Math. Soc. 358, no. 3, 1347-1364 (2006)

\bibitem[CCS2]{ccs2} Caldero P., Chapoton F., Schiffler
  R. \emph{Quivers with Relations and Cluster Tilted Algebras} Algebr
  Represent Theor 9, 359-376 (2006)

\bibitem[CK1]{ck1} Caldero P., Keller B. 
\emph{From triangulated categories to cluster algebras},
to appear in Invent. Math., arxiv: math.RT/0506018 

\bibitem[CK2]{ck2} Caldero P., Keller B. \emph{From triangulated
    categories to cluster algebras II}, Ann.
  Sci. \' Ecole Norm. Sup. (4) 39, no. 6, 983-1009 (2006)

\bibitem[CPS]{cps} Cline E., Parshall B., Scott L. \emph{Derived
    categories and Morita theory}, J. Algebra 104, no. 2, 397-409 (1986)

\bibitem[FZ1]{fz} Fomin S., Zelevinsky A. \emph{Cluster Algebras I:
    Foundations}, J. Amer. Math. Soc. 15, no. 2, 497-529 (2002)

\bibitem[FZ2]{fz2} Fomin S., Zelevinsky A. \emph{Cluster Algebras II:
    Finite type classification}, Invent. Math. 154, no. 1, 63--121 (2003)

\bibitem[Ha]{happel} Happel D. \emph{Triangulated Categories in the
    Representation Theory of Finite Dimensional Algebras}, London
  Mathematical Society Lecture Note Series, 119, Cambridge University
  Press, (1988)

\bibitem[Ho]{holm} Holm T. \emph{Cartan determinants for gentle
    algebras}, Arch. Math. 85, 233-239 (2005)

\bibitem[K]{keller} Keller B. \emph{On triangulated orbit
    categories}, Documenta Math. 10, 551-581 (2005)

\bibitem[KR]{kr} Keller B., Reiten I. \emph{Cluster-tilted algebras
    are Gorenstein and stably Calabi-Yau}, To appear in Adv. Math.,
  preprint arxiv:math.RT/0512471

\bibitem[M]{murphy} Murphy G. PhD thesis in preparation

\bibitem[Rin]{ringel} Ringel C.M. 
\emph{Some Remarks Concerning Tilting Modules and Tilted Algebras. 
Origin. Relevance. Future. (An appendix to the Handbook of Tilting Theory)},
Cambridge University Press (2007), LMS Lecture Notes Series 332

\bibitem[S]{seven} Seven A. \emph{Recognizing cluster algebras of
    finite type}, Electron. J. Combin. 14, no. 1, Research Paper 3, 35
  pp. (electronic) (2007)

\bibitem[SkW]{skw} Skowro\' nski A., Waschb\" usch
  J. \emph{Representation-finite biserial algebras}, J. Reine
  Angew. Math. 345, 172-181 (1983)

\end{thebibliography}
\end{document}